\documentclass[11pt]{amsart}
\usepackage[dvipsnames]{xcolor}
\definecolor{cornellred}{rgb}{0.7, 0.11, 0.11}
\newcommand{\redd}{\textcolor{cornellred}}
\newcommand{\greenn}{\textcolor{ForestGreen}}
\newcommand{\bluee}
{\textcolor{RoyalBlue}}
\usepackage[colorlinks=true,urlcolor=RoyalBlue,
citecolor=cornellred,linkcolor=RoyalBlue,linktocpage,pdfpagelabels,
bookmarksnumbered,bookmarksopen]{hyperref}

\usepackage{amssymb}
\usepackage{amsmath}
\usepackage{stmaryrd}
\usepackage{amsfonts}
\usepackage{enumerate}
\usepackage[shortlabels]{enumitem}
\setlist[enumerate]{align=left, leftmargin=*, labelsep=1 ex, parsep=.5ex, topsep=1ex}
\usepackage{subcaption}
\usepackage{comment}
\usepackage{tikz}
\usepackage{tikz-cd}

\definecolor{cornellred}{rgb}{0.7, 0.11, 0.11}
\newcommand{\red}{\textcolor{cornellred}}
\newcommand{\green}{\textcolor{ForestGreen}}
\newcommand{\blue}
{\textcolor{RoyalBlue}}

\usepackage{orcidlink}
\tikzset{multimap/.tip={Glyph[glyph math command=multimap]},}
\usetikzlibrary{decorations.pathmorphing}
\usepackage{mathtools}
\mathtoolsset{showonlyrefs}
\usepackage{comment}

\usepackage{amsthm}
\theoremstyle{plain}
\newtheorem{Theorem}{Theorem}[section]
\newtheorem{lemma}[Theorem]{Lemma}
\newtheorem{proposition}[Theorem]{Proposition}
\newtheorem{Fact}[Theorem]{Fact}
\newtheorem{corollary}[Theorem]{Corollary}

\theoremstyle{definition}
\newtheorem{definition}[Theorem]{Definition}
\newtheorem{Example}[Theorem]{Example}

\theoremstyle{remark}
\newtheorem{Remark}[Theorem]{Remark}

\usepackage[all]{xy}

\newcommand{\N}{\mathbb{N}}
\newcommand{\Q}{\mathbb{Q}}
\newcommand{\R}{\mathbb{R}}
\newcommand{\Cx}{\mathbb{C}}
\DeclareMathAlphabet{\mathmybb}{U}{bbold}{m}{n}

\DeclareMathOperator{\Ima}{Im}
\DeclareMathOperator{\Id}{id}
\DeclareMathOperator{\Ob}{Ob}
\DeclareMathOperator{\Sub}{\mathtt{Sub}}
\DeclareMathOperator{\Rel}{\mathtt{Rel}}
\DeclareMathOperator{\glb}{glb}
\newcommand{\Mag}{\mathtt{Mag}}
\DeclareMathOperator{\lub}{lub}
\DeclareMathOperator{\Map}{\mathtt{Map}}
\DeclareMathOperator{\Mon}{\mathtt{Mon}}
\DeclareMathOperator{\BMon}{\mathtt{BMon}}
\DeclareMathOperator{\Comon}{\mathtt{Comon}}
\DeclareMathOperator{\CMon}{\mathtt{CMon}}
\DeclareMathOperator{\Nwedge}{Nak_{\wedge}}
\DeclareMathOperator{\At}{At}
\DeclareMathOperator{\CComon}{\mathtt{CComon}}
\DeclareMathOperator{\Aut}{Aut}
\DeclareMathOperator{\Nvee}{Nak_{\vee}}
\DeclareMathOperator{\Hop}{\mathtt{Hop}}

\newcommand{\C}{\mathcal{C}}
\newcommand{\K}{\mathcal{K}}
\newcommand{\Msc}{\mathtt{Msc}}
\newcommand{\CPol}{\mathtt{CPol}}
\newcommand{\Mscs}{\mathtt{Msc}_\textup{str}}
\newcommand{\CPols}{\mathtt{CPol}_\textup{str}}
\newcommand{\LMsc}{\mathtt{LMsc}}
\newcommand{\CMsc}{\mathtt{CMsc}}
\newcommand{\CMscs}{\mathtt{CMsc}_\textup{str}}
\newcommand{\uMag}{\mathtt{uMag}}
\newcommand{\uMags}{\mathtt{uMag}_\textup{str}}
\newcommand{\Set}{\mathtt{Set}}
\newcommand{\Grp}{\mathtt{Grp}}
\newcommand{\Ab}{\mathtt{Ab}}

\title{L-mosaics and orthomodular lattices}
\author{Nicolò Cangiotti \orcidlink{0000-0003-2200-446X}}
\address[N. Cangiotti]{Dipartimento di Matematica, Politecnico di Milano, via Bonardi 9, Campus Leonardo, 20133 Milan, Italy}
\email{nicolo.cangiotti@polimi.it}
\author{Alessandro Linzi \orcidlink{0000-0002-1422-7674}}
\address[A. Linzi]{Center for Information Technologies and Applied Mathematics, University of Nova Gorica, Slovenia}
\email{linzi.alessandro@gmail.com}
\author{Enrico Talotti}
\address[E. Talotti]{Center for Information Technologies and Applied Mathematics, University of Nova Gorica, Slovenia}
\email{enrico.talotti@ung.si}

\begin{document}

\begin{abstract}
In this paper, we introduce a class of hypercompositional structures called dualizable L-mosaics. We prove that their category is equivalent to that formed by ortholattices and we formulate an algebraic property characterizing orthomodularity, suggesting possible applications to quantum logic.
\end{abstract}

\maketitle

\section{Introduction}

We introduce the notion of \emph{L-mosaic}, by extending a point of view initially pushed forward by T.\ Nakano in his studies on modular lattices \cite{Nak67} and by allowing non-associative multivalued operations as introduced in a recent  work of S.\ Nakamura and M.\ L.\ Reyes \cite{NR23}. We find a natural equivalence of categories between the category formed by ortholattices and the category of the multivalued algebraic structures we call \emph{dualizable L-mosaics}. This study is a follow-up on a recent article by A.\ Jen\v cov\'a and G.\ Jen\v ca suggesting applications of multivalued algebraic structures in quantum logic \cite{JJ17}.

The paper is organized as follows. In Sections \ref{Sec2} and \ref{Sec3}, we introduce mosaics, L-mosaics and ortholattices, while we collect a number of fundamental results on their structure theory and review the important property of ortholattices called orthomodularity. In Section \ref{Sec4} we study a particular class of L-mosaic which we call \emph{Nakano mosaics}. These L-mosaics are associated to any bounded lattice and are dualizable when the lattice is an ortholattice. In Section \ref{Sec5} we prove the equivalence of the categories mentioned above and formule the property of Nakano mosaics corresponding to orthomodularity. The connections with quantum logic and, more generally, with quantum mechanics are suggested in the conclusive Section \ref{Sec6}.

\section{Mosaics}
\label{Sec2}

Following Nakamura and M. L. Reyes \cite{NR23}, we introduce a multivalued algebraic structure called mosaic and review in detail some of their basic properties. 

\medskip

We consider two sets $A$ and $B$. By a \emph{multimap} $A\multimap B$ we mean a function $A\to\wp(B)$, with domain $A$ and codomain the power set $\wp(B)$ of $B$.

We call a multimap $f:A\multimap B$
\begin{enumerate}
\item \emph{total} if $f(x)\neq\emptyset$, for all $x\in A$,
\item \emph{partial} if it is not total,
\item \emph{deterministic} if $\lvert f(x)\rvert\leq 1$, for all $x\in A$,
\item a \emph{map} if it is total and deterministic.
\end{enumerate}

The composition $g\circ f:A\multimap C$ of two multimaps $f:A\multimap B$ and $g:B\multimap C$ is defined on every $x\in A$ by the following formula:
\[
(g\circ f)(x):=\bigcup_{y\in f(x)}g(y).
\]
It is easily verified that this composition is associative and that the multimaps defined by the assignment $x\mapsto\{x\}$ serve as identities. Moreover, it is not difficult to verify that the obtained category formed by sets and multimaps as above is isomorphic to the familiar category $\Rel$ formed by sets and binary relations. The isomorphism is provided by the equivalence
\begin{equation}\label{isomorphismR}
y\in f(x)\iff (x,y)\in R,
\end{equation}
where $f:A\multimap B$ is a multimap and $R\subseteq A\times B$ denotes the binary relation corresponding to $f$. By identifying these categories, we say that a multimap $f:A\multimap B$ \emph{contains} a multimap $g:A\multimap B$ if the relation corresponding to $f$ contains the relation corresponding to $g$.

A \emph{(binary) multioperation} on a set $A$ is a multimap ${\boxdot:A^2:=A\times A\multimap A}$. By a \emph{magma} we mean a set $A$ equipped with a multioperation.

To any multioperation $\boxdot$ defined on a set $A$ there correspond a multioperation \emph{dual} of $\boxdot$, defined for $x,y\in A$ by the formula
\[
x\boxdot^d y:=y\boxdot x.
\]

We call a magma $(A,\boxdot)$ \emph{total}/\emph{partial}/\emph{deterministic} if so is the multimap $\boxdot$. The magma $(A,\boxdot)$ is \emph{commutative} if the multioperations $\boxdot$ and $\boxdot^d$ coincide~on~$A$. If $(A,\boxdot)$ is total and deterministic, then we call it \emph{classical}.

By a \emph{morphism (of magmata)} from $(A,\boxdot_A)$ to $(B,\boxdot_B)$ we mean a function $f:A\to B$ such that, for all $x,y\in A$, the following inclusion holds:
\begin{equation}\label{hom}
f(x\boxdot_A y)\subseteq f(x)\boxdot_B f(y).
\end{equation}
By a \emph{strong morphism (of magmata)} from $(A,\boxdot_A)$ to $(B,\boxdot_B)$ we mean a function $f:A\to B$ such that, for all $x,y\in A$, the following inclusion holds:
\begin{equation}\label{shom}
f(x\boxdot_A y)= f(x)\boxdot_B f(y).
\end{equation}
An injective morphism of magmata from $f:(A,\boxdot_A)\to(B,\boxdot_B)$ is called an \emph{embedding} if 
\begin{equation}\label{emb}
f(x\boxdot_A y)=(f(x)\boxdot_B f(y))\cap f(A).
\end{equation}
It is easily seen that there we have a category $\Mag$ formed by magmata and morphisms and that isomorphisms in $\Mag$ are precisely bijective strong morphisms or, equivalently, surjective embeddings.

\begin{definition}
An element $e$ in a magma $(A,\boxdot)$ is called \emph{neutral} if $e\boxdot x=x\boxdot e=\{x\}$ holds, for all $x\in A$. A magma $(A,\boxdot)$ with a neutral element $e\in A$ for $\boxdot$ is called \emph{unital magma}.
\end{definition}

\begin{Remark}
If a neutral element $e$ exists in a magma $(A,\boxdot)$, then it is unique.
\end{Remark}

\begin{definition}
A morphism $f:A\to A'$ between two unital magmata $A$ and $A'$ with neutral elements $e$ and $e'$, respectively, is called \emph{unitary} if $f(e)=e'$.
\end{definition}

\begin{definition}
Let $(A,\boxdot)$ be a magma and $\rho:A\to A$ an endofunction. We say that $(A,\boxdot)$ is \emph{reversible with respect to $\rho$} (briefly, \emph{$\rho$-reversible}) if
\begin{enumerate}[label=\textup{(RE)}]
\item\label{rev} $z\in x\boxdot y$ implies both $x\in z\boxdot \rho(y)$ and $y\in \rho(x)\boxdot z$, for all $x,y,z\in A$.
\end{enumerate}
\end{definition}

\begin{definition}[{\cite[Definition 2.3]{NR23}}]
Let $(A,\boxdot,e)$ be a (commutative) unital magma,  which is $\rho$-reversible with respect to some endofunction $\rho:A\to A$ is called \emph{(commutative) mosaic}. By a \emph{(strong) morphism of mosaics} we mean a unitary (strong) morphism of the underlying unital magmata.
\end{definition}

The following is also observed in \cite{NR23}. We make a slightly more precise statement and write down a short proof.

\begin{lemma}[\cite{NR23}]\label{uniquerho}
Let $(A,\boxdot,e,\rho)$ be a mosaic. Then $\rho$ is a unitary isomorphism of magmata
\begin{equation}\label{isotrho}
\rho:(A,\boxdot,e)\overset{\sim}{\to}(A,\boxdot^d,e)
\end{equation}
satisfying the following property:
\begin{enumerate}[label=\textup{(RINV)}]
\item $e\in (x\boxdot \rho(x))\cap(\rho(x)\boxdot x)$, for all $x\in A$.\label{reinv}
\end{enumerate}
In addition, the equivalences
\[
z\in x\boxdot y\iff x\in z\boxdot\rho(y)\iff y\in \rho(x)\boxdot z  
\]
hold for all $x,y,z\in A$.
\end{lemma}

\begin{proof}
Indeed, $x\in (x\boxdot e)\cap(e\boxdot x)$ and $\rho$-reversibility imply $e\in (\rho(x)\boxdot x)\cap(x\boxdot \rho(x))$. Conversely, if $e\in (y\boxdot x)\cap(x\boxdot y)$ for some $y\in A$, then by $\rho$-reversibility we may deduce that $y\in e\boxdot \rho(x)=\{\rho(x)\}$ and hence $y=\rho(x)$. It follows that $\rho(e)=e$ and that $\rho$ is an involution. In addition, the validity of the following equivalences is readily verified, for all $x,y,z\in A$:
\begin{align*}
\rho(z)\in\rho(x)\boxdot^d \rho(y)&\iff \rho(z)\in\rho(y)\boxdot\rho(x)\\
&\iff \rho(y)\in \rho(z)\boxdot \rho(\rho(x))=\rho(z)\boxdot x\\
&\iff x\in \rho(\rho(z))\boxdot\rho(y)=z\boxdot \rho(y)\\
&\iff z\in x\boxdot \rho(\rho(y))=x\boxdot y.
\end{align*}
This shows that $\rho$ is a strong morphism of magmata $(A,\boxdot)\to (A,\boxdot^d)$ and thus an isomorphism. The rest of the assertions follow as well.
\end{proof}

\begin{definition}
For an element $x$ in a unital magma $(A,\boxdot,e)$, we call any $y\in A$ such that ${e\in (x\boxdot y)\cap(y\boxdot x)}$ an \emph{inverse} of $x$. If in $(A,\boxdot,e)$ all elements have a unique inverse, then we denote by $x^{-1}$ the inverse of any $x\in A$ and call $(A,\boxdot,e)$ an \emph{invertible magma}.
\end{definition}

The following is an immediate consequence of Lemma \ref{uniquerho}.

\begin{corollary}\label{invmag}
Let $(A,\boxdot,e,\rho)$ be a mosaic. Then $\rho(x)$ is the unique inverse of $x$ in $A$, for all $x\in A$. 
\end{corollary}

\begin{lemma}\label{to-1}
Let $(A,\boxdot,e)$ and $(A',\boxdot',e')$ be invertible magmata and $f:A\to A'$ a unitary morphism. Then $f(x^{-1})=f(x)^{-1}$, for all $x\in A$.
\end{lemma}

\begin{proof}
The standard argument as, e.g., in \cite[Remark 2.6]{KLS22}, applies.
\end{proof}

\begin{definition}[\cite{Com84}]
A \emph{(commutative) polygroup} $(P,\boxdot,e)$ is an invertible (commutative) magma, which is reversible with respect to the endofunction defined as ${x\mapsto x^{-1}}$ and where $\boxdot$ is \emph{associative}, that is, the following property is valid:
\begin{enumerate}[label=\textup{(ASC)}]
\item\label{asc} $(x\boxdot y)\boxdot z=x\boxdot (y\boxdot z)$, for all $x,y,z\in P$ 
\end{enumerate}
\end{definition}

\begin{Remark}
It follows from Corollary \ref{invmag} that polygroups are precisely associative mosaics.
\end{Remark}

\begin{definition}
By a \emph{(strong) morphism of polygroups} we mean a (strong) morphism of the underlying mosaics.
\end{definition}

The following result is another consequence of Lemma \ref{uniquerho}.

\begin{corollary}\label{rhorev}
Let $(A,\boxdot,e)$ be a unital magma satisfying \ref{asc}
and ${\rho:A\to A}$ an endofunction. The following statements are equivalent:
\begin{enumerate}
\item $(A,\boxdot,e)$ is $\rho$-reversible, i.e., $(A,\boxdot,e,\rho)$ is a polygroup.\label{co1}
\item\label{co2} $(A,\boxdot,e)$ is an invertible unital magma and $\rho(x)=x^{-1}$, for all $x\in A$.
\end{enumerate}
\end{corollary}

\begin{definition}
Let $(A,\boxdot,e)$ be a mosaic. A subset $B\subseteq A$ is a \emph{submosaic} if the inclusion map $B\to A$ is an embedding. A submosaic $B$ of a mosaic $(A,\boxdot,e)$ is \emph{strong} if the inclusion map $B\to A$ is a strong embedding.
\end{definition}

Examples of polygroups and mosaics can be found e.g.\ in \cite{BB19,NR23}.

\begin{definition}
A commutative mosaic $(A,\boxplus,0)$ is called a \emph{L-mosaic} if
\begin{enumerate}[label=\textup{(Lms\arabic*)}]
    \item $0,x\in x\boxplus x$, for all $x\in A$.\label{lm1}
    \item $(x\boxplus x)\boxplus (x\boxplus x)=x\boxplus x$, for all $x\in A$.\label{lm2}
    \item $(x\boxplus(x\boxplus y))\cap((x\boxplus y)\boxplus y)\subseteq x\boxplus y$, for all $x,y\in A$.\label{lm3}
    \item For all $x,y\in A$ there is a unique $z\in x\boxplus y$ such that $x,y\in z\boxplus z$.\label{lm4}
\end{enumerate}
\end{definition}

\begin{lemma}
    Let $(A,\boxplus,0)$ be an L-mosaic and $B$ a submosaic of $A$. Then the mosaic $(B,\boxplus_B,0)$ satisfies \ref{lm1}, \ref{lm2} and \ref{lm3}. 
\end{lemma}

\begin{proof}
    Obvious.
\end{proof}

\begin{definition}
    Let $(A,\boxplus,0)$ be an L-mosaic. A submosaic $B$ of $A$ is called \emph{L-submosaic} if $(B,\boxplus_B,0)$ satisfies property \ref{lm4}.
\end{definition}

\begin{lemma}\label{Lposm}
    Let $(A,\boxplus,0)$ be an L-mosaic. Then 
    \[
y\leq x\overset{\textup{\scriptsize{def}}}{\iff} y\in x\boxplus x.
\]
is an order relation on $A$ with respect to which $0$ is a bottom element.
\end{lemma}

\begin{proof}
    Reflexivity and the assertion on $0$ can be deduced immediately from \ref{lm1}. Regarding transitivity, note that if $x\in y\boxplus y$ and $y\in z\boxplus z$, then by \ref{lm2} implies 
    \[
    x\in y\boxplus y\subseteq (z\boxplus z)\boxplus (z\boxplus z)=z\boxplus z.
    \]
    As for antisymmetry, if $y\in x\boxplus x$ and $x\in y\boxplus y$, then by reversibility we have $x,y\in x\boxplus y$. Thus, property \ref{lm4} now implies $x=y$. 
\end{proof}

\begin{lemma}\label{Ax}
    Let $(A,\boxplus,0)$ be a mosaic satisfying properties \ref{lm1} and \ref{lm2}. If $x\in A$, then $A_x:=x\boxplus x$ is a strong submosaic of $A$. In addition, for any strong submosaic $B$ of $A$ we have that
    \[
    B=\bigcup_{x\in B} A_x=\bigcup_{x\in B} B_x.
    \]
\end{lemma}

\begin{proof}
    Let $x\in A$, then by \ref{lm1} we deduce that $0\in A_x$. By property \ref{lm2}, for all $y,z\in A_x$ we obtain that
    \[
    y\boxplus z\subseteq (x\boxplus x)\boxplus(x\boxplus x)=A_x.
    \]
    This suffices to prove that $A_x$ is a strong submosaic of $A$. In addition, we obtain that $B_x=A_x$, for all strong submosaics $B$ of $A$, whence
    \[
    \bigcup_{x\in B} A_x=\bigcup_{x\in B} B_x\subseteq B
    \]
 For the converse inclusion, note that by property \ref{lm1} it follows that $x\in B_x$ for all $x\in B$.
\end{proof}

The following is an obvious corollary.

\begin{corollary}\label{strcl}
    Let $(A,\boxplus,0)$ be a mosaic satisfying properties \ref{lm1} and \ref{lm2} and $B$ a submosaic of $A$. Then 
    \[
    \bar{B}:=\bigcup_{x\in B} A_x
    \]
    is a strong submosaic of $A$, which is contained in any strong submosaic of $A$ containing $B$.
\end{corollary}

\begin{definition}
    We call $\bar{B}$ the \emph{strong closure} of a submosaic $B$ in a mosaic $(A,\boxplus,0)$ satisfying properties \ref{lm1} and \ref{lm2}.
\end{definition}

\begin{definition}
    A commutative mosaic $(A,\boxplus,0)$ will be called \emph{$\pi$-dualizable} if $\pi:A\to A$ is an involution such that $(A,\boxplus_\pi,\pi(0))$ is a commutative mosaic, where for $x,y\in A$ we have set
    \[
    x\boxplus_\pi y:=\pi(x)\boxplus\pi(y).
    \]
\end{definition}

Clearly, $\pi$ becomes an isomorphism between the mosaic $(A,\boxplus,0)$ and its \emph{$\pi$-dual} $(A,\boxplus_\pi,\pi(0))$, whenever $A$ is $\pi$-dualizable. 

Therefore, a $\pi$-dualizable mosaic is an L-mosaic if and only if its $\pi$-dual is.

\section{Orthomodular lattices}
\label{Sec3}
We now review the main concepts around the orthomodularity property in the framework of lattices. Despite their introduction in the 1930s is strictly linked with the arising of quantum mechanics and the studies of J.\ V.\ Neumann \cite{VN32}, orthomodular lattices have become an independent purely algebraic dimension beyond their physical significance. We provide definitions in order to fix notations and make this paper as self-consistent as possible. The reader interested in further details is referred to  \cite{kalmbach83}.

A set $L$ equipped with two binary operations $\wedge$ and $\vee$ is a \emph{lattice} if
\begin{enumerate}[label=\textup{(L\arabic*)}]
\item Both $\wedge$ and $\vee$ are associative and commutative.
\item The equalities $x\vee(x\wedge y)=x$ and $x\wedge (x\vee y)=x$ hold, for all $x,y\in L$.
\end{enumerate}
A \emph{morphism of lattices $(L,\wedge,\vee)$, $(L',\wedge'\vee')$} is a function $f:L\to L'$ preserving operations. An \emph{isomorphism} is a bijective morphism.
A subset $L'$ of a lattice $L$ is called a \emph{sublattice} if it is closed under $\wedge$ and $\vee$. A \emph{bottom element} in a lattice $L$ is an element $0\in L$ such that $0 \leq x$ holds, for all $x\in L$. A lattice with a bottom element is called \emph{bounded from below}.  
A \emph{top element} in a lattice $L$ is an element $1\in L$ such that $x\leq 1$ holds, for all $x\in L$. A lattice with a top element is called \emph{bounded from above}.
A lattice with both a top and a bottom element is called \emph{bounded}.
Two elements $x,y$ in a bounded lattice $(L,\wedge,\vee,0,1)$ are \emph{complements}, written $x\bowtie y$, if and only if
\[
x\vee y=1\quad\textup{and}\quad x\wedge y=0.
\] 

Let $(L,\wedge,\vee)$ be a lattice. It can be seen {\cite[Chpt. 1, Sect. 2]{kalmbach83}} that the following relation
\[
x\leq y\overset{\textup{\scriptsize{def}}}{\iff} x\wedge y=x\iff x\vee y=y
\]
is a partial order on $L$ (called the \emph{canonical partial order} of the lattice $L$). In addition, for all $x,y,z\in L$ the following implications hold:
\begin{enumerate}[label=\textit{(\roman*)}]
\item\label{compi} $x\leq y\implies x\wedge z\leq y\wedge z$
\item\label{comps} $x\leq y\implies x\vee z\leq y\vee z$
\end{enumerate}

\begin{Remark}
Lattices resemble the algebraic properties of partially ordered sets where any pair of elements has an infimum and a supremum.  A priori, the roles of $\wedge$ and $\vee$ in a lattice are interchangeable. However, when the canonical order is defined a choice is made for the direction of the order, so that $\wedge$ denotes infima while $\vee$ denotes suprema.  The lattice obtained from a lattice $L$ by exchanging the roles of $\wedge$ and $\vee$, i.e., whose canonical order is the reverse order of the canonical order of $L$ will be referred to as the lattice \emph{dual to $L$}, denoted by $L^*$.

In addition, any theorem on lattices yields a \emph{dual theorem} obtained by exchanging the roles of $\vee$ and $\wedge$.
\end{Remark}

Let $(L,\wedge,\vee,0,1)$ be a bounded lattice. For $x\in L$, we define the multimap $\omega:L\multimap L$ by setting
\[
\omega(x):=\{y\in L\mid x\bowtie y\}.
\]
By the symmetry of the relation $\bowtie$, we obtain the following properties: 
\begin{enumerate}[label=\textup{(O\arabic*)}]
\item\label{omega1} $\omega(x)\neq\emptyset\implies x\in (\omega\circ\omega)(x)=:\omega^2(x)$. 
\item\label{omega2} $\omega^3(x)=\omega(x)$, for all $x\in L$.
\end{enumerate}
\begin{definition}
A bounded lattice $L$ is called \emph{complemented} if $\omega$ is total. A lattice $L$ is called an \emph{ortholattice} if $\omega$ contains a map (called \emph{orthocomplementation}).
\end{definition}

It follows from properties \ref{omega1} and \ref{omega2} above that all orthocomplementations on a bounded lattice are involutions.

\begin{Fact}[{\cite[Chpt. 1, Sect. 2]{kalmbach83}}]
Let $L$ be an ortholattice and $\pi:L\to L$ an orthocomplementation. Then the following statements hold:
\begin{enumerate}[label=\textup{(OC\arabic*)}]
\item $x\leq y$ if and only if $\pi(y)\leq \pi(x)$, for all $x,y\in L$.
\item $\pi(x\vee y)=\pi(x)\wedge\pi(y)$, for all $x,y\in L$.
\item $\pi(x\wedge y)=\pi(x)\vee\pi(y)$, for all $x,y\in L$.
\end{enumerate}
\end{Fact}

\begin{definition}
Let $(L,\wedge,\vee,0,1)$ be an ortholattice. If $\pi:L\to L$ is an orthocomplementation of $L$ such that for all $x,y\in L$ we have that
\[
x\leq y\implies x\vee (\pi(x)\wedge y)=y.
\]
Then the lattice $L$ is called \emph{orthomodular with respect to $\pi$} (\emph{$\pi$-orthomodular}). 
\end{definition}
\begin{definition}
Two elements $x,y$ in a $\pi$-orthomodular lattice are called \emph{orthogonal}, written $x\perp y$, if $x\leq \pi(y)$. 
\end{definition}

\begin{Example}
The ortholattice $H=\{0,a,a',b,b',1\}$ in Figure \ref{Benz} is not orthomodular.
Indeed $b\leq a$ but
\[
b\vee(b'\wedge a)=b\vee 0=b\neq a.
\]
\end{Example}

\begin{definition}
A lattice $(L,\wedge,\vee)$ is called \emph{modular} if for all $x,y,z\in L$ 
\[
x\vee(y\wedge (x\vee z))=(x\vee y)\wedge (x\vee z).
\]
\end{definition}

\begin{Fact}[{\cite[Chpt. 1, Sect. 2]{kalmbach83}}]\label{modulr}
A lattice $L$ is modular if and only if the pentagon lattice (Figure \ref{Pent}) is not a sublattice of $L$.
\end{Fact}

\begin{figure}[ht!]
\centering
    \begin{subfigure}[b]{.425\textwidth}
    \centering
    \begin{tikzpicture}[
C/.style = {circle, draw, fill=white, inner sep=1mm,
            label=#1,
            node contents={}},
every label/.append style = {inner sep=2pt}
                        ]
\draw (-2,1) node  (c1)     [C=left:$a$]    -- (0,2) 
             node   [C=$1$]         -- (2,1) node [C=right:$b'$] -- (2,-1) node [C=right:$a'$] -- (0,-2) node (c0) [C=below:$0$]
              ;
\draw (c1)  --  (-2,-1) node [C=left:$b$]
            --  (c0);
            \draw (c1)  --  (-2,-1) node [C=left:$b$]
            --  (c0);
    \end{tikzpicture}
    \caption{The hexagon lattice $H$.}
    \label{Benz}
    \end{subfigure}
    \qquad\qquad
    \begin{subfigure}[b]{.425\textwidth}
\begin{tikzpicture}[
C/.style = {circle, draw, fill=white, inner sep=1mm,
            label=#1,
            node contents={}},
every label/.append style = {inner sep=2pt}
                        ]
\draw (-2,1) node  (c1)     [C=left:$a$]    -- (0,2) 
             node   [C=$1$]         -- (2,0)
             node       [C=right:$c$]   -- (0,-2)
             node (c0)  [C=below:$0$]     -- (0,-2)
             ;
\draw (c1)  --  (-2,-1) node [C=left:$b$]
            --  (c0);
    \end{tikzpicture}
    \caption{The pentagon lattice $P$.}
    \label{Pent}
    \end{subfigure}
    
\end{figure}

\begin{Example}
The ortholattice $H=\{0,a,a',b,b',1\}$ in Figure \ref{Benz} is not orthomodular.
Indeed $b\leq a$ but
\[
b\vee(b'\wedge a)=b\vee 0=b\neq a.
\]
\end{Example}

The following characterization of orthomodularity is well-known.

\begin{Fact}[\cite{kalmbach83}, Theorem 2]\label{om}
\label{orthomchar}
    Let $L$ be a $\pi$-ortholattice. The following statements are equivalent:
    \begin{enumerate}[label=\textup{(OM\arabic*)}]
    \item $L$ is orthomodular.
    \item If $x\leq y$ and $\pi(x)\vee y=1$, then $x=y$, for all $x,y\in L$.\label{orthomcharii}
    \item\label{om3} $H$ is not a sublattice of $L$.
    \item If $x\leq y$, then the smallest sublattice $L'$ of $L$ containing $x,y$ and closed under $\pi$ is distributive.
    \end{enumerate}
\end{Fact}

Since the pentagon lattice $P$ of Figure \ref{Pent} is a sublattice of the hexagon lattice $H$ of Figure \ref{Benz}, it follows from Fact \ref{modulr} and Fact \ref{om} \ref{om3} that modularity is a property stronger than orthomodularity. For an example of an orthomodular lattice which is not modular \cite{kalmbach83}.

\begin{Fact}[\cite{kalmbach83}, Theorem 9]
\label{suborthmod}
Let $L$ be a $\pi$-orthomodular lattice and $x,y\in L$. Then the interstection of all sublattices $L'$ of $L$ such that $x,y\in L'$ and $\pi(L')=L'$ is a modular lattice.
\end{Fact}
\section{Nakano mosaics}
\label{Sec4}

We introduce in this section the objects we shall call \emph{Nakano mosaics}. These are inspired by the work of Nakano \cite{Nak67}, where the modularity of a lattice $L$ is shown to be equivalent to the associativity of some multioperations defined~on~$L$.

Throughout the section, a bounded lattice $(L,\wedge,\vee,0,1)$ is fixed and for all $x,y\in L$ we define the following subsets of $L$:
\begin{align*}
\Nvee(x,y)&:=\{z\in L\mid x\vee y=x\vee z=z\vee y\}\\
\Nwedge(x,y)&:=\{z\in L\mid x\wedge y=x\wedge z=z\wedge y\}
\end{align*}

\begin{proposition}\label{nakmos}
For all $x,y\in L$ set $x\boxplus y:=\Nvee(x,y)$. Then $(L,\boxplus,0)$ is a commutative, total and reproductive mosaic, where the inverse of each $x\in L$ is $x$ itself. Moreover, the top element $1$ can be characterized as the unique element $u\in L$ such that for all $x\in L$ the following implication holds
\[
u\in x\boxplus x\implies x=u.
\]
\end{proposition}

\begin{proof}
By definition $\Nvee(x,0)=\Nvee(0,x)=\{x\}$ holds for all $x\in L$. Moreover, it is easily verified from the definitions that 
\[
0\in \Nvee(x,y)\implies y=x,
\]
for all $x,y\in L$. Thus, $(L,\boxplus,0)$ is an invertible magma. To show that it is a mosaic it suffices to note that for $x,y,z\in L$ we have that:
\[
z\in x\boxplus y\iff x\vee y=z\vee x=z\vee y\iff x\in y\boxplus z.
\]
It is reproductive because, for fixed $x\in L$ and any $y\in L$, we have that, if $y\leq x$, then $y\in x\boxplus x$, while if $y>x$, then $y\in x\boxplus y$. In both cases, we obtain that for all $y\in L$ there is $x'\in L$ such that $y\in x\boxplus x'$, i.e., $L=x\boxplus L$. The assertion on the top element follows easily.
\end{proof}

\begin{proposition}[Dual of Proposition \ref{nakmos}]
For all $x,y\in L$ set $x\boxdot y:=\Nwedge(x,y)$. Then $(L,\boxdot,1)$ is a commutative, total and reproductive mosaic, where the inverse of each $x\in L$ is $x$ itself. Moreover, the bottom element $0$ can be characterized as the unique element $z\in L$ such that for all $x\in L$
\[
z\in x\boxdot x\implies x=z.
\]
\end{proposition}

\begin{proposition}\label{subl+}
Let $L$ be a bounded lattice. The additive Nakano mosaic associated to a sublattice $L'$ of $L$ such that $0\in L'$ is a submosaic of the additive Nakano mosaic associated to $L$.
\end{proposition}

\begin{proof}
Let $L'$ be a sublattice of $L$ such that $0\in L'$. Consider the additive Nakano mosaic $(L',\boxplus',0)$ associated to $L'$. We show that $(x \boxplus' y)=(x\boxplus y)\cap L'$, for all $x,y\in L$. The assertion clearly follows. In particular, we have $z\in x\boxplus y$ if and only if $x\vee y=x\vee z=y\vee z$, hence $z\in x\boxplus'y$ if and only if $z\in (x\boxplus y)\cap L'$.
\end{proof}

\begin{proposition}[Dual of Proposition \ref{subl+}]\label{subl.}
Let $L$ be a bounded lattice. The multiplicative Nakano mosaic associated to a sublattice $L'$ of $L$ such that $1\in L'$ is a submosaic of the multiplicative Nakano mosaic associated to $L$.
\end{proposition}

\begin{lemma}\label{antisym}
Let $(L,\boxplus,0)$ be the additive Nakano mosaic associated to the bounded lattice $L$. The following equivalence
\[
x,y\in x\boxplus y\iff x=y.
\]
holds, for all $x,y\in L$
\end{lemma}

\begin{proof}
We have $x,y\in x\boxplus y$ if and only if both $x\in y\boxplus y$ and $y\in x\boxplus x$ hold. The assertion thus follows by the antisymmetry of the canonical order $\leq$ of the lattice $L$.
\end{proof}

\begin{lemma}\label{trans}
Let $(L,\boxplus,0)$ be the additive Nakano mosaic associated to the bounded lattice $L$, then
\[
(x\boxplus x)\boxplus (x\boxplus x)= x\boxplus x
\]
holds, for all $x\in L$.
\end{lemma}

\begin{proof}
If $y,z\leq x$ holds in $L$, then $t\leq t\vee z=t\vee y=z\vee y\leq x$ follows, for all $t\in y\boxplus z$, i.e., $t\in x\boxplus x$. The converse inclusion is immediately verified after noticing that $x\in x\boxplus x$ holds for all $x\in L$.
\end{proof}

\begin{Fact}[{\cite[Lemma 24]{Corsini03}}]
\label{CL}
Let $(L,\boxplus,0)$ be the additive Nakano mosaic associated to the bounded lattice $L$ and $x,y,z\in L$. Then
\[
x\boxplus (y\boxplus z)\subseteq\{t\in L\mid t\vee x\vee y=t\vee x\vee z=t\vee y\vee z=x\vee y\vee z\}.
\]
\end{Fact}

\begin{corollary}\label{lemlm3}
Let $(L,\boxplus,0)$ be the additive Nakano mosaic associated to the bounded lattice $L$, then
\[
(x\boxplus (x\boxplus y))\cap ((x\boxplus y)\boxplus y)\subseteq x\boxplus y
\]
holds, for all $x,y\in L$.
\end{corollary}

\begin{proof}
Follows by Fact \ref{CL} and the properties of the lattice operation $\vee$.
\end{proof}

\begin{lemma}\label{leqsup}
For all $x,y,z\in L$, we have $z\in x\boxplus y$ implies $z\leq x\vee y$. Thus, $x\boxplus y\subseteq (x\vee y)\boxplus (x\vee y)$. 
\end{lemma}

\begin{proof}
    If $z\in x\boxplus y$, then
    \[
    z\vee(x\vee y)=(z\vee x)\vee y=(x\vee y)\vee y=x\vee(y\vee y)=x\vee y,
    \]
    hence $z\leq x\vee y$.
\end{proof}

\begin{lemma}\label{sup}
For $x,y,z\in L$ we have that $z=x\vee y$ if and only if
\begin{equation}\label{eqqr}
x,y\in z\boxplus z\textup{ and }z\in x\boxplus y
\end{equation}
hold in the additive Nakano mosaic $(L,\boxplus,0)$ associated to the bounded lattice $L$.
\end{lemma}

\begin{proof}
If $z=x\vee y$, then $x,y\leq z$ and consequently $x\vee y=z=x\vee z=y\vee z$, which shows \eqref{eqqr}. Conversely, if \eqref{eqqr} holds, then (e.g.) $x\leq z$ follows and thus $x\vee y=x\vee z=z$.   
\end{proof}

\begin{corollary}\label{+Lmos}
    The additive Nakano mosaic associated to a bounded lattice is an $L$-mosaic.
\end{corollary}

\begin{lemma}[Dual of Lemma \ref{leqsup}]\label{geqinf}
For all $x,y,z\in L$, we have $z\in x\boxdot y$ implies $z\geq x\wedge y$. Thus, $x\boxdot y\subseteq (x\wedge y)\boxdot (x\wedge y)$. 
\end{lemma}

\begin{lemma}[Dual of Lemma \ref{sup}]\label{inf}
For $x,y,z\in L$ we have that $z=x\wedge y$ if and only if
\[
x,y\in z\boxdot z\textup{ and }z\in x\boxdot y
\]
hold in the multiplicative Nakano mosaic $(L,\boxdot,1)$ associated to the bounded lattice $L$.
\end{lemma}

\begin{corollary}
    The multiplicative Nakano mosaic associated to a bounded lattice is an $L$-mosaic.
\end{corollary}

\begin{Fact}[{\cite[Theorem 1]{Nak67}}]
\label{originalNak}
For a bounded lattice $L$, the following are equivalent statements:
\begin{enumerate}
    \item $(L,\boxplus,0)$ is a polygroup.
    \item $(L,\boxdot,1)$ is a polygroup.
    \item $L$ is a modular lattice.
\end{enumerate}
\end{Fact}

\begin{Example}
The Pentagon lattice $P=\{0,a,b,c,1\}$ in Figure \ref{Pent} is not modular. Hence, the associated Nakano mosaics $\mathbf{P}_\vee$ and $\mathbf{P}_\wedge$ are not associative. Table \ref{Pentvtab} and Table \ref{Pentwtab} fully describe the multioperations $\boxplus$ of $\mathbf{P}_\vee$ and $\boxdot$ of $\mathbf{P}_\wedge$.  Note that $\mathbf{P}_\vee$ is not a polygroup as $a\boxplus (b\boxplus c)=\{c,1\}\neq\{1\}=(a\boxplus b)\boxplus c$. Moreover $\mathbf{P}_\vee$ is not isomorphic to $\mathbf{P}_\wedge$.
\begin{table}[ht]
    \centering
    \begin{tabular}{|c|c|c|c|c|c|}
    \hline
      $\boxplus$  & $0$ &$a$ & $b$ & $c$ & $1$  \\ \hline
       $0$ & $\{0\}$ & $\{a\}$ & $\{b\}$&$\{c\}$ & $\{1\}$  \\ \hline
       $a$ & $\{a\}$ & $\{0,a,b\}$ & $\{a\}$ & $\{1\}$& $\{c,1\}$ \\ \hline
       $b$ & $\{b\}$& $\{a\}$ & $\{0,b\}$& $\{1\}$ & $\{c,1\}$ \\ \hline
       $c$ &$\{c\}$ & $\{1\}$ & $\{1\}$ &$\{0,c\}$ & $\{a,b,1\}$\\ \hline
       $1$ & $\{1\}$& $\{c,1\}$ & $\{c,1\}$ & $\{a,b,1\}$ & $\{0,a,b,c,1\}$ \\      \hline 
    \end{tabular}
    \caption{$\mathbf{P}_\vee$: The additive Nakano mosaic associated to the non-modular Pentagon lattice.}
    \label{Pentvtab}
    \end{table}
    \qquad\qquad
\begin{table}[ht]
\centering
    \begin{tabular}{|c|c|c|c|c|c|}
    \hline
      $\boxdot$  & $1$ &$a$ & $b$ & $c$ & $0$  \\ \hline
       $1$ & $\{1\}$ & $\{a\}$ & $\{b\}$&$\{c\}$ & $\{0\}$  \\ \hline
       $a$ & $\{1,a\}$ & $\{1\}$ & $\{b\}$ & $\{0\}$& $\{c,0\}$ \\ \hline
       $b$ & $\{b\}$& $\{a,b\}$ & $\{1,b\}$& $\{0\}$ & $\{c,0\}$ \\ \hline
       $c$ &$\{c\}$ & $\{0\}$ & $\{0\}$ &$\{1,c\}$ & $\{a,b,0\}$\\ \hline
       $0$ & $\{0\}$& $\{c,0\}$ & $\{c,0\}$ & $\{a,b,0\}$ & $\{0,a,b,c,0\}$ \\      \hline 
    \end{tabular}
    \caption{$\mathbf{P}_\wedge$: The multiplicative Nakano mosaic associated to the non-modular Pentagon lattice.}
    \label{Pentwtab}
    \end{table}
\end{Example}

\begin{Example}
    In Table \ref{Benztab} we describe the multioperation $\boxplus$ of the additive Nakano mosaic associated to the hexagon lattice $H$ of Figure \ref{Benz}.
    
\begin{table}[ht]
    \centering
    \begin{tabular}{|c|c|c|c|c|c|c|}
    \hline
      $\boxplus$  & $0$ &$a$ & $b$ & $a'$ & $b'$ & $1$  \\ \hline
       $0$ & $\{0\}$ & $\{a\}$ & $\{b\}$&$\{a'\}$ & $\{b'\}$ & $\{1\}$  \\ \hline
       $a$ & $\{a\}$ & $\{0,a,b\}$ & $\{a\}$ & $\{1\}$ & $\{1\}$& $\{a',b',1\}$ \\ \hline
       $b$ & $\{b\}$& $\{a\}$ & $\{0,b\}$& $\{1\}$ & $\{1\}$& $\{a',b',1\}$ \\ \hline
       $a'$ &$\{a'\}$ & $\{1\}$ & $\{1\}$ &$\{0,a'\}$ & $\{b'\}$& $\{a,b,1\}$\\ \hline
       $b'$ &$\{b'\}$ & $\{1\}$ & $\{1\}$ &$\{a'\}$ & $\{0,a',b'\}$& $\{a,b,1\}$\\ \hline
       $1$ & $\{1\}$& $\{b,c\}$ & $\{a,c\}$ & $\{a,b\}$ & $\{b'\}$& $\{0,a,b,c,1\}$ \\      \hline 
    \end{tabular}
    \caption{$\mathbf{H}_\vee$: The additive Nakano polygroup associated to the hexagon ortholattice.}
    \label{Benztab}
\end{table}
\end{Example}

\begin{corollary}
Let $L$ be a bounded lattice. If $\mathbf{P}_\vee$ is not a submosaic of $(L,\boxplus,0)$, then the additive Nakano mosaic $(L,\boxplus,0)$ associated to $L$ is a polygroup.
\end{corollary}

\begin{proof}
If $(L,\boxplus,0)$ is not a polygroup, then $L$ is not modular by Fact \ref{originalNak}, thus the Pentagon lattice is a sublattice of it. By Proposition \ref{subl+}, it follows that $\mathbf{P}_\vee$ is a submosaic of $(L,\boxplus,0)$.   
\end{proof}

\begin{Remark}
    Clearly, by composition with the duality endofunctor on the category of bounded lattices, we obtain another fully faithful functor: mapping a bounded lattice to its multiplicative Nakano mosaic.
\end{Remark}

Since fully faithful functors are conservative, we obtain a useful corollary.

\begin{corollary}\label{cons}
Let $(L,\wedge,\vee,0,1),(L',\wedge',\vee',0',1')$ be bounded lattices. The associated additive Nakano mosaics $(L,\boxplus,0)$ and $(L',\boxplus',0')$ are isomorphic if and only if the lattices $L$ and $L'$ are. 

By duality, both the associated multiplicative Nakano mosaics $(L,\boxdot,1)$ and $(L',\boxdot',1')$ are isomorphic if and only if the lattices $L$ and $L'$ are.      
\end{corollary}

\section{Orthomodularity and Nakano mosaics}
\label{Sec5}

Let $\mathtt{Ort}$ be the category formed by pairs $(L,\pi)$, where $L$ is a $\pi$-ortholattice, as objects and lattice morphisms $f:L\to L'$ such that $f\circ\pi=\pi'\circ f$ as arrows $f:(L,\pi)\to(L',\pi')$.

Let $\LMsc^d$ the category formed by pairs $(A,\pi)$, where $A$ is a $\pi$-duaizable L-mosaic as objects and mosaic morphisms $f:A\to A'$ such that $f\circ\pi=\pi'\circ f$ as arrows $f:(A,\pi)\to(A',\pi')$.

\begin{Theorem}\label{!!}
    The categories $\mathtt{Ort}$ and $\LMsc^d$ are equivalent.
\end{Theorem}

\begin{proof}

    We employ Proposition \ref{nakmos} to define the object assignment of a functor $\mathcal{E}$ as 
    \[
    \mathcal{E}(L,\wedge,\vee,0,1):=(L,\boxplus,0).
    \]

    We start by claiming that $f:L\to L'$ is a morphism of lattices if and only if the same map is a morphism of the corresponding mosaics. Indeed, if $z\in x\boxplus y$ holds in the mosaic $(L,\boxplus,0)$ defined above, then we deduce
    \begin{align*}
    f(x)\vee f(y)&=f(x\vee y)=f(x\vee z)=f(x)\vee f(z),\textup{ and}\\
    f(x)\vee f(y)&=f(x\vee y)=f(z\vee y)=f(z)\vee f(y),
    \end{align*}
    that is, $f(z)\in f(x)\boxplus f(y)$. 

    Conversely, if $f$ is a morphism of mosaics, then since $x\vee y\in x\boxplus y$ holds by definition of $\boxplus$ for any $x,y\in L$, we obtain that $f(x\vee y)\in f(x)\boxplus' f(y)$, whence $f(x\vee y)\vee f(x)=f(x)\vee f(y)$. On the other hand, we have that $x\vee y\geq x$, i.e., $x\vee y\in x\boxplus x$. Therefore, $f(x\vee y)\in f(x)\boxplus' f(x)$, meaning that $f(x\vee y)\geq f(x)$ holds in $L'$. We conclude that $f(x\vee y)\vee f(x)=f(x\vee y)=f(x)\vee f(y)$. The equality $f(x\wedge y)=f(x)\wedge f(y)$ is obtained employing the duality isomorphism $\pi$ and the fact that $x\wedge y=\pi(\pi(x)\vee\pi(y))$. 

    Thus, we have a fully faithful functor $\mathcal{E}:\mathtt{Ort}\to\LMsc^d$, it remains to prove that it is essentially surjective, i.e., if $(A,\boxplus,0)$ is a $\pi$-dualizable L-mosaic, then we have to define a unique (up to isomorphism) lattice structure on $A$ such that $\pi$ is an orthocomplementation and $(A,\boxplus,0)$ is the associated (additive) Nakano mosaic. Since $(A,\boxplus,0)$ is an L-mosaic, by Lemma \ref{Lposm} we obtain the poset $(A,\leq)$ with bottom element~$0$.

    We claim that for all $x,y\in A$ the unique $z$, whose existence is guaranteed by property \ref{lm4}, is the least upper bound of $x$ and $y$ in $(A,\leq)$. Indeed, it is an upper bound since $x,y\in z\boxplus z$. If $z'\in A$ is any upper bound, i.e, $x,y\in z'\boxplus z'$, then $z\leq z'$ because
    \[
    z\in x\boxplus y\subseteq (z'\boxplus z')\boxplus (z'\boxplus z')=z'\boxplus z'.
    \]
    Now, using the assumption that $(A,\boxplus,0)$ is $\pi$-dualizable we set $x\vee y:=z$ and $x\wedge y:=\pi(\pi(x)\vee \pi(y))$ for all $x,y\in A$ and we obtain a $\pi$-ortholattice $(A,\wedge,\vee,0,\pi(0))$. It remains to verify that $a\boxplus b=\Nvee(a,b)$ for all $a,b$ in the lattice $(A,\wedge,\vee,0,\pi(0))$.

    If $c\in a\boxplus b$ holds in $(A,\boxplus,0)$ and we set $x:=a\vee b$, $y:=a\vee c$ as well as $z:=b\vee c$, then we have that $a\in x\boxplus x$ and 
    \[
    c\in a\boxplus b\subseteq (x\boxplus x)\boxplus (x\boxplus x)=x\boxplus x.
    \]
    Thus, we deduce $a,c\leq x$ and $y=a\vee c\leq x$. On the other hand,
    \[
    b\in a\boxplus c\subseteq (y\boxplus y)\boxplus(y\boxplus y)=y\boxplus y.
    \]
    Thus, we deduce $a,b\leq y$ and $x=a\vee b\leq y$. Therefore, $x=y$ follows by antisymmetry and similarly one proves that $z=x=y$. This proves that $a\boxplus b\subseteq \Nvee(a,b)$.

    For the converse inclusion, take $c\in \Nvee(a,b)$, i.e., $a\vee b=a\vee c=b\vee c=:x$. By definition of $\vee$, we obtain that 
    \[
    x\in (a\boxplus b)\cap(a\boxplus c)\cap(b\boxplus c),
    \]
    therefore, using reversibility,
    \[
    c\in a\boxplus x\subseteq a\boxplus (a\boxplus b)
    \]
    and 
    \[
    c\in x\boxplus b\subseteq (a\boxplus b)\boxplus b
    \]
    $c\in a\boxplus b$ follows by property \ref{lm3} of L-mosaics.
    \end{proof}

\begin{proposition}
Let $A$ be a $\pi$-dualizable L-mosaic. Then the $\pi$-ortholattice associated to $A$ is $\pi$-orthomodular if and only if the implication 
\begin{equation}\label{omodLmos}
x\in y\boxplus y\textup{ and }1\in x\boxplus \pi(y)\implies x=y
\end{equation}
holds, for all $x,y\in A$. 
\end{proposition}

\begin{proof}
    We recall that, by Fact \ref{orthomchar} \ref{orthomcharii}, the $\pi$-ortholattice associated to $A$ is $\pi$-orthomodular if and only if the conjunction of $x\leq y$ and $x\vee\pi(y)=1$ implies $x=y$. By definition of $\vee$ in the $\pi$-ortholattice associated to $A$ and since $1\boxplus 1=A$, we have that $x\vee\pi(y)=1$ is equivalent to $1\in x\boxplus \pi(y)$. On the other hand, we know that $x\leq y$ means that $x\in y\boxplus y$. The assertion follows.
\end{proof}

\begin{proposition}
    Let $A$ be a $\pi$-dualizable L-mosaic satisfying the implication \eqref{omodLmos} for all $x,y\in A$. Then for all $x,y\in A$ the smallest L-submosaic $A'$ of $A$ containing $x,y$ and closed under $\pi$ is a polygroup.
\end{proposition}
\begin{proof}
    This is a straightforward consequence of Fact \ref{suborthmod}.
\end{proof}

\section{Further research: the connection with the quantum world}
\label{Sec6}

Since the 1990s, \emph{dagger compact categories} have played a central role in the development of topological quantum field theories, as introduced by John Baez and James Dolan \cite{Baez95}. Later, in the early 2000s, Samson Abramsky and Bob Coecke identified these categories as fundamental structures in their framework of categorical quantum mechanics \cite{AC04}.

Building on these foundational ideas, recent work has explored alternative algebraic approaches to quantum theories. In this context, we highlight the contribution of Jen\v{c}ov\'a and Jen\v{c}a \cite{JJ17}, which served as an initial inspiration for our study. Their research underscores the potential role of hypercompositional structures in the analysis of quantum phenomena. Specifically, their approach employs the concept of \emph{effect algebras}—algebraic structures with a partial operation (see \cite{FB94} for further details). 

Both partial operations and multioperations can be naturally modeled by monoids in the category $\Rel$ (see \cite{LinAIMS} for further details). At this point, it should be noted that the category $\Rel$ possess a dagger compact structure. While partial algebraic structures, such as effect algebras, are already well-integrated into quantum theories, the study of algebraic structures with multioperations remains less developed. This work aims to encourage further research into this promising direction, which we also plan to explore in future investigations.

\section*{Acknowledgments}
The authors are grateful towards the organizers and  participants of both the \emph{ ``4$^{th}$ Symposium On Hypercompositional Algebra-New Developments And Applications''} and  the \emph{\lq\lq 20th International Conference on Quantum Physics and Logic\rq\rq}, where the research contained in this paper has been fruitfully discussed.


\begin{thebibliography}{99}

  
\bibitem{AC04}
S.~Abramsky and B.~Coecke., A categorical semantics of quantum protocols, 
  \emph{Proceedings of the 19th Annual IEEE Symposium on Logic in Computer Science},  415--425, 2004.

\bibitem{Baez95}
J.~C. Baez and J.~Dolan, Higher-dimensional Algebra and Topological Quantum Field Theory, \emph{J. Math. Phys.}, 36(11): 6073–-6105, 1995. 

\bibitem{BB19}
M.~Baker and N.~Bowler.
\newblock Matroids over partial hyperstructures.
\newblock {\em Adv. Math.}, 343:821--863, 2019.


\bibitem{Com84}
S.~D. Comer.
\newblock Combinatorial aspects of relations.
\newblock {\em Algebra Universalis}, 18(1):77--94, 1984.


\bibitem{Corsini03}
P.~Corsini and V.~Leorenau.
\newblock {\em Applications of Hyperstructure Theory}.
\newblock Springer, New York, USA, 2003.
 

\bibitem{FB94}
D.~J. Foulis and M.~K. Bennett.
\newblock Effect algebras and unsharp quantum logics.
\newblock {\em Found. Phys.}, 24:1331--1352, 1994.


\bibitem{JJ17}
A.~Jen\v{c}ov\'{a} and G.~Jen\v{c}a.
\newblock On monoids in the category of sets and relations.
\newblock {\em Internat. J. Theoret. Phys.}, 56(12):3757--3769, 2017.


\bibitem{kalmbach83}
G.~Kalmbach.
\newblock {\em Orthomodular Lattices}.
\newblock L.M.S. monographs. Academic Press, 1983.

\bibitem{KLS22}
K.~Kuhlmann, A.~Linzi, and H.~Stoja{\l}owska.
\newblock Orderings and valuations in hyperfields.
\newblock {\em J. Algebra}, 611:399--421, 2022.

\bibitem{LinAIMS}
A.~Linzi.
\newblock Polygroup objects in regular categories.
\newblock {\em AIMS Mathematics}, 9(5):11247-11277, 2024.


\bibitem{NR23}
S.~{N}akamura and M.~L. {R}eyes.
\newblock Categories of hypermagmas, hypergroups, and related hyperstructures.
\newblock {a}r{X}iv:2304.09273, 2023.

\bibitem{Nak67}
T.~Nakano.
\newblock Rings and partly ordered systems.
\newblock {\em Math. Z.}, 99:355--376, 1967.


\bibitem{VN32}
J. Von Neumann.
\newblock {\em Mathematical Foundations of Quantum Mechanics}.
\newblock Vol. 183 Princeton University Press, Princeton, 1932.
\end{thebibliography}
\end{document}